 \documentclass[11pt]{article}
 \usepackage{indentfirst, latexsym,bm}
 \usepackage{amsmath}
 \usepackage{pifont}
 \usepackage{amsfonts}
 \usepackage{mathrsfs}
 \usepackage{array}
 \usepackage{multirow}
 \usepackage{graphicx}
 \usepackage{subfigure}
 \usepackage{picinpar}

 \usepackage{booktabs}
 \usepackage{threeparttable}
 \usepackage{mathrsfs,amsfonts,amsmath}
 \usepackage[dvips]{color}

\newcommand\email[1]{\href{mailto:#1}{ \nolinkurl{#1}}}

 \newtheorem{theorem}{Theorem}[section]
 \newtheorem{definition}[theorem]{Definition}
 \newtheorem{lemma}[theorem]{Lemma}
 \newtheorem{corollary}[theorem]{Corollary}
 \newtheorem{proposition}[theorem]{Proposition}
 \newtheorem{remark}[theorem]{Remark}
 \newtheorem{condition}[theorem]{Condition}
 \newtheorem{example}{Example}[section]

 \def\blemma{\begin{lemma}}\def\elemma{\end{lemma}}
 \def\bproposition{\begin{proposition}}\def\eproposition{\end{proposition}}
 \def\ttheorem{\begin{theorem}}\def\etheorem{\end{theorem}}
 \def\bcorollary{\begin{corollary}}\def\ecorollary{\end{corollary}}
 \def\bremark{\begin{remark}}\def\eremark{\end{remark}}
 \def\bcondition{\begin{condition}}\def\econdition{\end{condition}}

 \def\benumerate{\begin{enumerate}}\def\eenumerate{\end{enumerate}}
 \def\bitemize{\begin{itemize}}\def\eitemize{\end{itemize}}
 \def\itm{\item}

 \def\beqlb{\begin{eqnarray}}\def\eeqlb{\end{eqnarray}}
 \def\beqnn{\begin{eqnarray*}}\def\eeqnn{\end{eqnarray*}}
 \def\nnm{\nonumber}\def\ccr{\nnm\\}\def\ar{\!\!\!&}

 \def\mcr{\mathscr}\def\mbb{\mathbb}\def\mbf{\mathbf}

 \def\proof{\noindent{\it Proof.~~}}\def\qed{\hfill$\Box$\medskip}

\begin{document}

\noindent{(Version: 2016/02/05)}

\bigskip\bigskip

\centerline{\LARGE\bf Asymptotic results for exponential}

\centerline{\LARGE\bf functionals of L\'{e}vy processes\,\footnote{Supported by the NSFC
(No.~11131003, No.~11531001  and  No.~11401012).}}

\bigskip

\centerline{Zenghu Li and Wei Xu}

\medskip

\centerline{School of Mathematical Sciences, Beijing Normal University,}

\centerline{Beijing 100875, People's Republic of China}

\centerline{E-mails: \tt lizh@bnu.edu.cn, xuwei@mail.bnu.edu.cn}

\bigskip

{\narrower{\narrower

\textit{Abstract.} The long time asymptotic behavior of the expectation of some exponential functional of a L\'{e}vy process is studied. We give not only the exact convergence rate but also explicitly the limiting coefficients. The key of the results is the observation that the asymptotics only depends on the sample paths of the L\'evy process with local infimum decreasing slowly. This makes it possible for us to determine the limiting coefficients by extending the conditional limit theorems for L\'evy processes established by Hirano (2001). The constants are represented in terms of some transformations based on the renewal functions. As applications of the results, we give exact evaluation of the decay rate of the survival probability of a continuous-state branching process in random environment with stable branching mechanism.

\medskip


\medskip


\par}\par}

\bigskip

\section{Introduction}

\setcounter{equation}{0}

The study of exponential functionals of random walks and L\'evy processes has drawn the attention of many researchers in recent years. Those functionals play important roles in the study of probabilistic models in random environments among their other applications. Let $\xi = \{\xi(t): t\ge 0\}$ be a one-dimensional L\'evy process. Given a constant $\alpha>0$, we define the exponential functional:
 \beqlb\label{eq1.1}
A_t^\alpha(\xi) = \int_0^t e^{-\alpha \xi(s)}ds, \qquad 0\le t\le \infty.
 \eeqlb
When $\xi$ is a Brownian motion with drift, a characterization of the distribution of $A_t^\alpha(\xi)$ was obtained by Yor (1992, Proposition~2). For an exponentially distributed random variable $T$, positive and negative moments of $A_T^\alpha(\xi)$ were calculated by Carmona et al.\ (1994, 1997). By a result of Bertoin and Yor (2005), we have $A_\infty^\alpha(\xi)< \infty$ if and only if $\lim_{t\to \infty} \xi(t) = \infty$. In this case, Bertoin and Yor (2005) gave some characterizations for the distribution of $A_\infty^\alpha(\xi)$ and Pardo et al.\ (2012) established a Wiener-Hopf type factorization for this functional. Let $z\mapsto F(z)$ be a positive decreasing function on $(0,\infty)$ that vanishes as $z\to \infty$ at a certain rate. In the case of $A_\infty^\alpha(\xi) =\infty$, a natural problem is to evaluate the decay rate as $t\to \infty$ of the expectation:
 \beqlb\label{eq1.2}
\mbf{P}[F(A_t^\alpha(\xi))]
 =
\mbf{P}\bigg[F\Big(\int_0^t e^{-\alpha \xi(s)}ds\Big)\bigg].
 \eeqlb
In the special case where $F(z) = a(a+z)^{-1}$ and $\{\xi(t): t\ge 0\}$ is a Brownian motion with drift, the problem was studied by Kawazu and Tanaka (1993) in their work on the tail behavior of a diffusion process in random environment. Other specific forms of the function $F$ arising from applications were discussed in Carmona et al.\ (1994, 1997).

Let $\{Z_\alpha(t): t\ge 0\}$ be a spectrally positive $(\alpha+1)$-stable processes with $0<\alpha\le 1$ and $\{L(t): t\ge 0\}$ a L\'{e}vy process with no jump less than $-1$. Let $c\ge 0$ be another constant. Given the initial value $x\ge 0$, we consider the following stochastic integral equation:
 \beqlb\label{eq1.3}
X(t) = x + \int_0^t \sqrt[1+\alpha]{(1+\alpha)cX(s-)} dZ_\alpha(s) + \int_0^t X(s-)dL(s).
 \eeqlb
By Theorem~6.2 in Fu and Li (2010), there exists a unique positive strong solution $\{X(t): t\ge 0\}$ to (\ref{eq1.3}). The solution is called a \textit{continuous-state branching process in random environment (CBRE-process)} with \textit{stable branching mechanism}. Here the random environment is modeled by the L\'evy process $\{L(t): t\ge 0\}$. The reader may refer to He et al.\ (2016) and Palau and Pardo (2015b) for discussions of more general CBRE-processes. We shall see that there is another L\'{e}vy process $\{\xi(t): t\ge 0\}$ determined by the environment so that the \textit{survival probability} of the CBRE-process up to time $t\ge 0$ is given by
 \beqlb\label{eq1.4}
\mathbf{P}(X(t)>0) = \mathbf{P}\Big[1-\exp\big\{-x(c\alpha)^{-1/\alpha}A_t^\alpha(\xi)^{-1/\alpha}\big\}\Big].
 \eeqlb
Clearly, the right-hand side of (\ref{eq1.4}) is a special case of (\ref{eq1.2}). Based on the above expression, the asymptotic behavior of the survival probability as $t\to \infty$ was studied by B\"{o}inghoff and Hutzenthaler (2012) for the case where $\alpha=1$ and the environment process is a Brownian motion with drift. Their results were extended recently to the case $1<\alpha\le 1$ by Palau and Pardo (2015a). The main strategy of B\"{o}inghoff and Hutzenthaler (2012) and Palau and Pardo (2015a) is the formula of Yor (1992) for the distribution of the exponential functional of the Brownian motion with drift; see also Matsumoto and Yor (2003). Bansaye et al.\ (2013) studied the problem in the case where the environment is given by a L\'evy process with bounded variations and showed some interesting applications of the results to a cell infection model. The key step in their proof is to study the expectation (\ref{eq1.2}) for $F(x) = (1+x)^{-1/\beta} [1+(1+x)^{-\gamma} h(x)]$, where $0<\beta\le 1$ and $\gamma\ge 1$ are constants and $h$ is a bounded Lipschitz function. The asymptotics of survival probabilities for classical Galton-Watson branching processes in random environment (GWRE-processes) was studied earlier by Afanasy'ev et al.\ (2005), Dyakonova et al.\ (2004), Geiger and Kersting (2002), Geiger et al.\ (2003), Guivarc'h and Liu (2001), Kozlov (1976), Liu (1996) and Vatutin et al.\ (2013) among others. Roughly speaking, for critical branching the survival probability decays at a polynomial rate and for subcritical branching it decays at an exponential rate with three different polynomial modifying factors, which classify the processes into weakly subcritical, intermediately subcritical and strongly subcritical ones. Those results play important roles in the study of various conditional limit theorems of the CBRE- and GWRE-processes. Unfortunately, in most of the results established before, the limiting coefficients were not explicitly identified except in very special cases; see, e.g., B\"{o}inghoff and Hutzenthaler (2012).

The purpose of this paper is to study the asymptotic behavior of the expectation in (\ref{eq1.2}) for a general function $F$ and a general L\'evy process $\xi$. Under natural assumptions, we prove some accurate results for asymptotics of the expectation as $t\to \infty$. We shall see that five regimes arise for the convergence rate. We also apply the results to study the survival probability of the CBRE-process defined by (\ref{eq1.3}). The feature of this work is that we give explicitly not only the convergence rate but also the limiting coefficient in all regimes. The key of the results is the observation that the asymptotics of (\ref{eq1.2}) only depends on the charge of probability on sample paths of the L\'evy process whose local infimum decreases slowly. This makes it possible for us to determine the limiting coefficients by extensions of the conditional limit theorems of Hirano (2001). The constants are represented in terms of some transformations based on the renewal functions associated with the ladder processes of $\xi$ and its dual process. The main results of the paper are presented in Section~2. The proofs for recurrent and transient L\'evy processes are given in Sections~3 and~4, respectively. The applications of the results to CBRE-processes are discussed in Section~5.

After putting the first version of this paper to Arxiv, we noticed the interesting work of Palau et al.\ (2016), where some results for the asymptotics of exponential functionals of L\'evy processes were obtained and the results were also applied to study the survival probability of the CBRE-process. But, as in most of the references mentioned above, they did not identify the limiting coefficients.

\medskip

\textbf{Acknowledgements.} We thank Professor Gennady Samorodnitsky for enlightening comments on the literature of L\'evy processes. We are grateful to Professors Juan C. Pardo and Vladimir A. Vatutin for letting us know their work on branching processes in random environments.

\section{Asymptotics of exponential functionals}\label{ExpLevy}

\setcounter{equation}{0}

In this section, we study the asymptotic behavior of some exponential functionals of L\'{e}vy processes. We first introduce some basic notations. Let ${\it\Psi}(\lambda)$ the \textit{characteristic exponent} of an infinitely divisible probability measure on $\mbb{R}$ given by
 \beqlb\label{eq2.1}
{\it\Psi}(\lambda) = ia\lambda + \frac{\sigma^2}{2}\lambda^2 + \int_{\mbb{R}}(1-e^{i\lambda x}+i\lambda x)\nu(dx), \qquad \lambda\in \mbb{R},
 \eeqlb
where $a\in \mbb{R}$ and $\sigma\ge 0$ are constants and $\nu(dx)$ is a $\sigma$-finite measure on $\mbb{R}$ supposed by $\mbb{R}\setminus \{0\}$ and satisfying
 \beqlb\label{eq2.2}
\int_{\mbb{R}}(|x|\wedge |x|^2)\nu(dx)<\infty.
 \eeqlb
We also need to consider the \textit{Laplace exponent} ${\it\Phi}$ defined by
 \beqlb\label{eq2.3}
{\it\Phi}(\lambda) = -a\lambda + \frac{\sigma^2}{2}\lambda^2 + \int_{\mbb{R}}(e^{\lambda x} - 1 - \lambda x)\nu(dx), \qquad \lambda\in \mbb{R}.
 \eeqlb
Of course, we may have ${\it\Phi}(\lambda) = \infty$ for some $\lambda\in \mbb{R}$. Let $\mcr{D}({\it\Phi}) = \{\lambda\in \mbb{R}: {\it\Phi}(\lambda)< \infty\}$. Then $\mcr{D}({\it\Phi})$ is necessarily an interval containing the origin. Let $\mcr{D}_+({\it\Phi}) = \mcr{D}({\it\Phi})\cap [0,\infty)$. Let $\mcr{D}^\circ({\it\Phi})$ and $\mcr{D}^\circ_+({\it\Phi})$ denote the interior sets of $\mcr{D}({\it\Phi})$ and $\mcr{D}_+({\it\Phi})$, respectively.

Let ${\it\Omega}$ be the set of all c\`adl\`ag paths from $[0,\infty)$ to $\mbb{R}$. For $t\ge 0$ and $\omega\in {\it\Omega}$ let $\xi_t(\omega) = \omega(t)$ denote the \textit{coordinate process}. Let $\mcr{F} = \sigma(\{\xi_t: t\ge 0\})$ and $\mcr{F}_t = \sigma(\{\xi_s: 0\le s\le t\})$ be the natural $\sigma$-algebras. For each $x\in \mbb{R}$ there is a probability measure $\mbf{P}_x$ on $({\it\Omega}, \mcr{F})$ so that $\{(\xi_t, \mcr{F}_t): t\ge 0\}$ under this measure is a process with independent and stationary increments and
 \beqlb\label{eq2.4}
\mbf{P}_x[\exp\{i\lambda \xi_t\}] = \exp\{i\lambda x-t{\it\Psi}(\lambda)\}, \qquad t\ge 0,\lambda\in \mbb{R}.
 \eeqlb
Then $\xi = ({\it\Omega}, \mcr{F}, \mcr{F}_t, \xi_t, \mbf{P}_x)$ is the \textit{canonical realization} of the L\'evy process with characteristic exponent ${\it\Psi}$. Let $\hat{\mbf{P}}_x$ denote the law of $\{-\xi_t: t\ge 0\}$ under $\mbf{P}_{-x}$. Then $\hat{\xi} = ({\it\Omega}, \mcr{F}, \mcr{F}_t, \xi_t, \hat{\mbf{P}}_x)$ is the \textit{dual process} of $\xi$, which is also a L\'evy process. For simplicity, write $\mbf{P} = \mbf{P}_0$ and $\hat{\mbf{P}} = \hat{\mbf{P}}_0$. It is well-known that under the integrability condition (\ref{eq2.2}) we have
 \beqnn
\mbf{P}_0(\xi_t) = -\hat{\mbf{P}}_0(\xi_t) = -at, \qquad t\ge 0.
 \eeqnn
For notational convenience, we also write $\xi(t)$ instead of $\xi_t$ in the sequel.

For any $\theta\in \mcr{D}({\it\Phi})$, it is easy to see that $t\mapsto e^{\theta \xi(t)-{\it\Phi}(\theta)t}$ is a $\mbf{P}_x$-martingale. Then, using the Escheer transform, we can define the probability measure $\mbf{P}_x^{(\theta)}$ on $({\it\Omega}, \mcr{F})$ by
 \beqlb\label{eq2.5}
\mbf{P}_x^{(\theta)}(A)
 =
\int_A e^{\theta \xi(t)-{\it\Phi}(\theta)t}d\mbf{P}_x, \qquad A\in \mcr{F}_t, t\ge 0.
 \eeqlb
It is known that $\xi^{(\theta)} = ({\it\Omega}, \mcr{F}, \mcr{F}_t, \xi(t), \mbf{P}_x^{(\theta)})$ is a L\'evy process with Laplace exponent ${\it\Phi}_\theta(\lambda) := {\it\Phi}(\lambda+\theta) - {\it\Phi}(\lambda)$; see, e.g., Theorem~3.9 in Kyprianou (2014, p.83). Let $\hat{\xi}^{(\theta)} = ({\it\Omega}, \mcr{F}, \mcr{F}_t, \xi(t), \hat{\mbf{P}}_x^{(\theta)})$ be its dual process.

We define the \textit{supremum process} $S := (S(t): t\ge 0)$ by $S(t) = \sup_{s\in[0,t]}\xi(s)$. Let $S-\xi := \{S(t)-\xi(t): t\ge 0\}$ be the \textit{reflected process}, which is a Markov process with Feller transition semigroup; see, e.g., Proposition~1 in Bertoin (1996, p.156). Let $L = \{L(t): t\ge 0\}$ be the local time at zero of $S-\xi$. The \textit{inverse local time process} $L^{-1} = \{L^{-1}(t): t\ge 0\}$ is defined by
 \beqnn
L^{-1}(t)=\Big\{
 \begin{array}{ll}\inf\{s>0: L(s)>t\}, & t< L(\infty);\cr
 \infty, & \mbox{otherwise}.
 \end{array}
 \eeqnn
The \textit{ladder height process} $H=\{H(t):t\ge 0\}$ of $\xi$ is defined by
 \beqnn
H(t)=\Big\{\begin{array}{ll}\xi(L^{-1}(t)), & t< L(\infty);\cr
 \infty, & \mbox{otherwise}.
 \end{array}
 \eeqnn
By Lemma~2 in Bertoin (1996, p.157), the two-dimensional process $(L^{-1},H)$ is a (possibly killed) L\'{e}vy process. This is known as the \textit{ladder process} of $\xi$ and is characterized by
 \beqnn
\mbf{P}[\exp\{-\lambda_1 L^{-1}(t) - \lambda_2 H(t)\}]
 =
\exp\{-t\kappa(\lambda_1,\lambda_2)\}, \qquad \lambda_1,\lambda_2\ge 0,
 \eeqnn
where the \textit{bivariate exponent} $\kappa(\lambda_1,\lambda_2)$ is given by
 \beqnn
\kappa(\lambda_1,\lambda_2)
 \ar=\ar
k\exp\Big\{\int_0^\infty \frac{dt}{t}\int_{[0,\infty)}(e^{-t}-e^{-\lambda_1 t-\lambda_2 x})\mbf{P}(\xi(t)\in dx)\Big\}.
 \eeqnn
The constant $k>0$ here is determined by the normalization of the local time; see Corollary~10 in Bertoin (1996, pp.165--166). In this work, we choose the normalization suitably so that $k=1$; see also Hirano (2001, p.293). The \textit{renewal function} $V$ associated with the ladder height process $H$ is defined by
 \beqlb\label{eq2.6}
V(x) = \int_0^\infty \mbf{P}(H(t)\le x)dt, \qquad x\ge 0.
 \eeqlb
Let $\hat{V}$ and $\hat{\kappa}(\lambda_1,\lambda_2)$ be defined similarly as the above from $\hat{\xi}$. Let $V(\cdot-)$ and $\hat{V}(\cdot-)$ denote the left limits of the renewal functions.

For $x\in \mbb{R}$ define the hitting time $\tau_x = \inf\{t> 0: \xi_t\le x\}$. It is known that for any $x>0$ the process $t\mapsto \hat{V}(\xi(t)-)\mbf{1}_{\{\tau_0> t\}}$ is a $\mbf{P}_x$-martingale and $t\mapsto V(\xi(t)-)\mbf{1}_{\{\tau_0> t\}}$ is a $\hat{\mbf{P}}_x$-martingales; see Bertoin (1996, p.184) and Hirano (2001, p.293). Then we can define the probability measures $\mbf{Q}_x$ and $\hat{\mbf{Q}}_x$ on $({\it\Omega}, \mcr{F})$ by
 \beqnn
\mbf{Q}_x(A)
 =
\hat{V}(x-)^{-1} \int_A\hat{V}(\xi(t)-)\mbf{1}_{\{\tau_0> t\}} d\mbf{P}_x
 \eeqnn
and
 \beqnn
\hat{\mbf{Q}}_x(A)
 =
V(x-)^{-1} \int_AV(\xi(t)-)\mbf{1}_{\{\tau_0> t\}} d\hat{\mbf{P}}_x,
 \eeqnn
where $A\in \mcr{F}_t$ and $t\ge 0$. Let ${\it\Xi} = ({\it\Omega}, \mcr{F}, \mcr{F}_t, \xi(t), \mbf{Q}_x)$ and $\hat{{\it\Xi}} = ({\it\Omega}, \mcr{F}, \mcr{F}_t, \xi(t), \hat{\mbf{Q}}_x)$.

Let $V^{(\theta)}$ and $\hat{V}^{(\theta)}$ be the renew functions associated with the ladder height processes of the L\'evy processes $\xi^{(\theta)}$ and $\hat{\xi}^{(\theta)}$, respectively. Let ${\it\Xi}^{(\theta)} = ({\it\Omega}, \mcr{F}, \mcr{F}_t, \xi(t), \mbf{Q}_x^{(\theta)})$ and $\hat{{\it\Xi}}^{(\theta)} = ({\it\Omega}, \mcr{F}, \mcr{F}_t, \xi(t), \hat{\mbf{Q}}_x^{(\theta)})$ be the resulted Markov processes, respectively.

For $\alpha>0$ let $\{A_t^\alpha(\xi): 0\le t\le \infty\}$ be defined by (\ref{eq1.1}).

\begin{proposition}\label{th2.1}{\rm (Carmona et al., 1997)} For any $\alpha>0$ the following statements are equivalent: (1) $\mbf{P}[\xi(1)]> 0$; (2) $\mbf{P}(A_\infty^\alpha(\xi)< \infty)>0$; (3) $\mbf{P}(A_\infty^\alpha(\xi)< \infty) = 1$. \end{proposition}

\begin{lemma}\label{th2.2} For any $\alpha>0$, $t>0$ and $\beta\in \mcr{D}_+^\circ({\it\Phi})$ we have
 \beqnn
\mbf{P}[A_t^\alpha(\xi)^{-\beta/\alpha}]
 \le
t^{-\beta/\alpha}\mbf{P}[e^{\beta S(t)}]
 \le
4t^{-\beta/\alpha}e^{\beta(a-|a|)t}\mbf{P}[e^{\beta \xi(t)}].
 \eeqnn
\end{lemma}

\proof By applying Doob's inequality to the submartingale $t\mapsto e^{\beta[\xi(t)+at]/2}$ we have
 \beqnn
\mbf{P}[e^{\beta S(t)}]
 \ar\le\ar
e^{\beta|a|t}\mbf{P}\Big[\sup_{0\le s\le t}e^{\beta[\xi(s)+as]}\Big]
 \le
e^{\beta|a|t}\mbf{P}\Big[\sup_{0\le s\le t}\Big(e^{\beta[\xi(s)+as]/2}\Big)^2\Big] \ccr
 \ar\le\ar
4e^{\beta|a|t}\mbf{P}[e^{\beta[\xi(t)+at]}]
 =
4e^{\beta(a+|a|)t}\mbf{P}[e^{\beta \xi(t)}],
 \eeqnn
where the right-hand side is finite since $\beta\in \mcr{D}_+^\circ({\it\Phi})$. It is simple to see that
 \beqnn
\mbf{P}[A_t^\alpha(\xi)^{-\beta/\alpha}]
 \le
\mbf{P}\Big[\Big(\int_0^t e^{-\alpha S(t)}ds\Big)^{-\beta/\alpha}\Big]
 =
t^{-\beta/\alpha}\mbf{P}[e^{\beta S(t)}]
 \eeqnn
Then we obtain the result. \qed

\begin{lemma}\label{th2.3} For any $\alpha>0$, $t>0$ and $\beta\in \mcr{D}_+^\circ({\it\Phi})$ we have
 \beqnn
\mbf{P}[A_t^\alpha(\xi)^{-\beta/\alpha}]
 \le
\mbf{P}\Big[\exp\Big\{\min_{k\le [t]-1}\beta \xi(k)\Big\}\Big]\mbf{P}[e^{\beta S(1)}].
 \eeqnn
\end{lemma}

\proof It suffices to consider $t\ge 2$ in this proof. Let $[t]$ denote the integer part of $t$. For any $j=0,1,\cdots, [t]-1$, define
 $$
Z(j)=\log\Big(\int_j^{j+1}e^{-\alpha(\xi(s)-\xi(j))}ds\Big).
 $$
Then $\{Z(j): j=0,1,\cdots,[t]-1\}$ is a sequence of i.i.d. random variables. It is easy to see that
 \beqnn
\mbf{P}[A_t^\alpha(\xi)^{-\beta/\alpha}]
 \ar\le\ar
\mbf{P}\Big[\Big(\int_0^{[t]} e^{-\alpha \xi(s)}ds\Big)^{-\beta/\alpha}\Big]\cr
 \ar=\ar
\mbf{P}\Big[\Big(\sum_{j=0}^{[t]-1} e^{-\alpha \xi(j)+Z(j)}\Big)^{-\beta/\alpha}\Big]
 \le
\mbf{P}[e^{\beta \xi(\kappa)-\beta Z(\kappa)/\alpha}],
 \eeqnn
where $\kappa=\min\{j\le [t]-1: \xi(j)=\min_{k\le [t]-1} \xi(k)\}$. Since $Z(\kappa)$ is independent of $\xi(\kappa)$ and $\kappa$, we have
 \beqnn
\mbf{P}[A_t^\alpha(\xi)^{-\beta/\alpha}]
 \ar\le\ar
\sum_{k=0}^{[t]-1} \mbf{P}(\kappa=k)\mbf{P}[e^{\beta \xi(\kappa)-\beta Z(\kappa)/\alpha}|\kappa=k]\cr
 \ar=\ar
\sum_{k=0}^{[t]-1} \mbf{P}(\kappa=k)\mbf{P}[e^{\beta \xi(\kappa)}|\kappa=k] \mbf{P}[e^{-\beta Z(\kappa)/\alpha}]\cr
 \ar=\ar
\sum_{k=0}^{[t]-1} \mbf{P}(\kappa=k)\mbf{P}[e^{\beta \xi(\kappa)}|\kappa=k] \mbf{P}[e^{-\beta Z(0)/\alpha}]\cr
 \ar\le\ar
\sum_{k=0}^{[t]-1} \mbf{P}(\kappa=k)\mbf{P}[e^{\beta \xi(\kappa)}|\kappa=k] \mbf{P}[e^{\beta S(1)}]\cr
 \ar=\ar
\mbf{P}\Big[\exp\Big\{\min_{k\le [t]-1}\beta \xi(k)\Big\}\Big] \mbf{P}[e^{\beta S(1)}].
 \eeqnn
Then the desired result follows. \qed

\begin{lemma}\label{th2.4} For any $\alpha>0$ and $x>0$ we have $\mbf{Q}_x[A_\infty^\alpha(\xi)]< \infty$. \end{lemma}

\proof By the definition of $\mbf{Q}_x$ and Fubini's theorem, we have
 \beqnn
\mbf{Q}_x[A_\infty^\alpha(\xi)]
 \ar=\ar
\hat{V}(x-)^{-1}\int_0^{\infty}\mbf{P}_x[e^{-\alpha \xi(r)}\hat{V}(\xi(r)-);\tau_0>r]dr \cr
 \ar=\ar
\hat{V}(x-)^{-1}\mbf{P}_x\Big[\int_0^{\tau_0}e^{-\alpha \xi(r)}\hat{V}(\xi(r)-)dr\Big] \cr
 \ar\le\ar
\hat{V}(x-)^{-1}\int_0^\infty dV(y)\int_0^x e^{-\alpha (y+x-z)}\hat{V}(y+x-z)d\hat{V}(z),
 \eeqnn
where the last step follows by Theorem~20 in Bertoin (1996, p.176). By Corollary~5.3 in Kyprianou (2014, p.118) we have $\hat{V}(y)\sim y/\hat{\mbf{P}}[H(1)]$ as $y\to \infty$. Then we can take $\gamma\in (0,\alpha)$ and $C\ge 0$ so that $e^{-(\alpha-\gamma)y} \hat{V}(y)\le C$ for $y\ge 0$. It follows that
 \beqnn
\mbf{Q}_x[A_\infty^\alpha(\xi)]
 \le
C\hat{V}(x-)^{-1}\int_0^\infty e^{-\gamma y} dV(y)\int_0^x d\hat{V}(z),
 \eeqnn
The right-hand side is clearly finite. \qed

 \medskip
\noindent{\textbf{General Assumption}\,} \textit{In the sequel of the paper, we fix constants $\alpha> 0$, $\beta\in \mcr{D}_+^\circ({\it\Phi})$ and a decreasing and strictly positive function $F$ on $(0,\infty)$. In addition, we assume there exist constants $C_0>0$ and $\beta_0\in \mcr{D}_+({\it\Phi})$ so that $F(z)\le C_0z^{-\beta_0/\alpha}$ for $0<z\le 1$.}
 \medskip

To simplify the presentation of the results, let us state the following conditions:

\begin{condition}\label{th2.5} For each $\delta>0$ there is a constant $K_\delta> 0$ so that $|F(z)-F(y)|\le K_\delta|z-y|$ for $z,y\ge \delta$. \end{condition}

\begin{condition}\label{th2.6} There is a constant $K>0$ so that $F(z)\le Kz^{-\beta/\alpha}$ for $z\ge 1$. \end{condition}

\begin{condition}\label{th2.7} There is a constant $K>0$ so that $F(z)\sim Kz^{-\beta/\alpha}$ as $z\to \infty$. \end{condition}

\begin{condition}\label{th2.8} The characteristic exponent satisfies $\mbf{Re}{\it\Psi}(\lambda)> 0$
for all $\lambda\neq 0$. \end{condition}

The main theorem of this paper is the following:

\begin{theorem}\label{th2.9}
 \begin{enumerate}

\itm[\rm(1)] If $\mbf{P}[\xi(1)]> 0$, we have the non-degenerate limit
 $$
\lim_{t\to \infty}\mbf{P}[F(A_t^\alpha(\xi))] = \mbf{P}[F(A_\infty^\alpha(\xi))].
 $$

\itm[\rm(2)] Suppose that Conditions~\ref{th2.5}, \ref{th2.6} and~\ref{th2.8} are satisfied. If $0\in \mcr{D}^\circ({\it\Phi})$ and ${\it\Phi}'(0)=0$, then we have the non-degenerate limit
 \beqnn
\lim_{t\to \infty}t^{1/2}\mbf{P}[F(A_t^\alpha(\xi))]
 =
\sqrt{\frac{2}{\pi{\it\Phi}''(0)}}\hat{\mbf{P}}[H(1)] D_2(\alpha,F),
 \eeqnn
where
 \beqlb\label{eq2.7}
D_2(\alpha,F) = \lim_{x\to \infty} \hat{V}(x-) \mbf{Q}_x[F(e^{-\alpha x}A_\infty^\alpha(\xi))].
 \eeqlb

\itm[\rm(3)] Suppose that Conditions~\ref{th2.5}, \ref{th2.6} and~\ref{th2.8} are satisfied and that $0\in \mcr{D}^\circ({\it\Phi})$ and ${\it\Phi}'(0)< 0< {\it\Phi}'(\beta)$. Let $\varrho\in (0,\beta)$ be the solution of ${\it\Phi}'(\varrho) = 0$ and let $(W, \mcr{G}, \mcr{G}_t, (\xi(t),\hat{\xi}(t)), \mbf{Q}^{(\varrho)}_{(x,y)})$ be the independent coupling of ${\it\Xi}^{(\varrho)}$ and $\hat{{\it\Xi}}^{(\varrho)}$. Then we have the non-degenerate limit
 \beqnn
\lim_{t\to\infty}t^{3/2}e^{-t{\it\Phi}(\varrho)}\mbf{P}[F(A_t^\alpha(\xi))]
 =
\frac{c(\varrho)}{\sqrt{2\pi{\it\Phi}''(\varrho)}} D_3(\alpha,F),
 \eeqnn
where
 \beqlb\label{eq2.8}
c(\varrho) = \exp\Big\{\int_0^\infty (e^{-t}-1)t^{-1}e^{-t{\it\Phi}(\varrho)} \mbf{P}(\xi(t)=0)dt\Big\},
 \eeqlb
 \beqlb\label{eq2.9}
D_3(\alpha,F) = \lim_{x\to \infty} e^{\varrho x} \hat{V}^{(\varrho)}(x-) \int_0^\infty e^{-\varrho y} V^{(\varrho)}(y) G(x,y)dy
 \eeqlb
and
\beqlb\label{eq2.10}
G(x,y) = \mbf{Q}^{(\varrho)}_{(x,y)} \{F(e^{-\alpha x}[A_\infty^\alpha(\xi)+A_\infty^\alpha(\hat{\xi})])\}.
 \eeqlb

\itm[\rm(4)] Suppose that Conditions~\ref{th2.7} and~\ref{th2.8} are satisfied and ${\it\Phi}'(\beta) = 0$. Then we have the non-degenerate limit
 \beqnn
\lim_{t\to \infty}t^{1/2}e^{-t{\it\Phi}(\beta)}\mbf{P}[F(A_t^\alpha(\xi))]
 =
K\sqrt{\frac{2}{\pi{\it\Phi}''(\beta)}} \mbf{P}^{(\beta)}[H(1)] D_4(\alpha,\beta),
 \eeqnn
where
 \beqlb\label{eq2.11}
D_4(\alpha,\beta) = \lim_{x\to \infty} V^{(\beta)}(x-) \mbf{Q}_x^{(\beta)} [e^{-\beta x}A_\infty^\alpha(-\xi)^{-\beta/\alpha}].
 \eeqlb

\itm[\rm(5)] Suppose that Condition~\ref{th2.7} is satisfied and ${\it\Phi}'(\beta)< 0$. Then we have the non-degenerate limit
 \beqnn
\lim_{t\to \infty}e^{-t{\it\Phi}(\beta)}\mbf{P}[F(A_t^\alpha(\xi))]
 =
K\mbf{P}^{(\beta)}[A_\infty^\alpha(-\xi)^{-\beta/\alpha}].
 \eeqnn

 \end{enumerate}
\end{theorem}

\begin{remark}\label{th2.10} It is known that Condition~\ref{th2.8} holds if and only if $\sigma> 0$ or $\nu(\mathbb{R}\backslash \{0, \pm r, \pm 2r, \cdots\})> 0$ for every $r>0$; see, e.g., Hirano (2001, p.294). Instead of this condition, if we assume for some $r>0$ the characteristic exponent has the representation:
 \beqlb\label{eq2.12}
{\it\Psi}(\lambda) = \sum_{k\in \mbb{Z}}(1-e^{ikr\lambda}), \qquad \lambda\in \mbb{R},
 \eeqlb
the results of regimes (2) and (4) in the above theorem still hold. The proofs are modifications of those given in Sections~3 and~4. However, it seems some extra work is needed to establish the result in regime (3) for the characteristic exponent (\ref{eq2.12}). \end{remark}

By using the above theorem we can give some simple derivations of the results of B\"{o}inghoff and Hutzenthaler (2012), Carmona et al.\ (1994, 1997) and Kawazu and Tanaka (1993) on the asymptotics of exponential functionals; see Xu (2016).

\section{Recurrent L\'evy processes}\label{Critical}

\setcounter{equation}{0}

In this section, we consider the case where the L\'evy process $\xi$ is recurrent. In particular, we shall give the proof of Theorem~\ref{th2.9} in regime (2). Throughout the section, we assume $0\in \mcr{D}^\circ({\it\Phi})$. It follows that $\mbf{P}[\xi(1)] = {\it\Phi}'(0) = 0$ and $\mbf{P}[\xi(1)^2] = {\it\Phi}''(0)< \infty$.

\begin{proposition}\label{th3.1}
 \begin{enumerate}

\itm[\rm(1)] Let $I(t) = \inf_{0\le s\le t} \xi(s)$. Then for $x>0$ and $\epsilon>0$ we have, as $t\to \infty$,
 \beqnn
\mbf{P}(\tau_{-x}>t) = \mbf{P}(I(t)>-x)\sim \sqrt{\frac{2}{\pi{\it\Phi}''(0)}}\hat{\mbf{P}}[H(1)]\hat{V}(x-) t^{-1/2}.
 \eeqnn
\itm[\rm(2)] Suppose that Condition~\ref{th2.8} holds. Then for any $x>0$ and $\alpha>0$ we have, as $t\to \infty$,
 \beqnn
\mbf{P}(e^{-\alpha \xi(t)}; \tau_{-x}>t)
 \sim
\frac{c e^{\alpha x}}{\sqrt{2\pi{\it\Phi}''(0)}}\hat{V}(x-) t^{-3/2}\int_{0}^\infty e^{-\alpha z}V(z)dz,
 \eeqnn
where
 \beqnn
c = \exp\Big\{\int_0^{\infty}(e^{-t}-1)t^{-1} \mbf{P}(\xi(t)=0)dt\Big\}.
 \eeqnn

 \end{enumerate}
\end{proposition}

\proof Since $\mbf{P}(I(t)>-x) = \hat{\mbf{P}}(\sup_{0\le s\le t}\xi(s)< x)$, the first result follows from Lemma~11 of Hirano (2001). By the spatial homogeneity of the L\'evy process we have
 \beqnn
\mbf{P}(e^{-\alpha \xi(t)};\tau_{-x}>t)
 =
\mbf{P}(e^{-\alpha[x+\xi(t)]}; \tau_{-x}>t)e^{\alpha x}
 =
\mbf{P}_x(e^{-\alpha \xi(t)}; \tau_0>t)e^{\alpha x}.
 \eeqnn
Then the second result follows by Lemma~A-(a) in Hirano (2001). \qed

For $s\ge 0$ let $D[0,s]$ denote the space of c\`adl\`ag real functions on $[0,s]$ equipped with the Skorokhod topology. The following proposition extends slightly Theorem~1 in Hirano (2001), who considered the case where $\xi(t)\to -\infty$ as $t\to \infty$. Based on Proposition~\ref{th3.1}(1), its proof goes similarly as that given in Hirano (2001), so we omit the proof here. The discrete version of the result for random walks was established by Bertoin and Doney (1994).

\begin{proposition}\label{th3.2} Let $f$ be a bounded Borel function on $D[0,s]$. Then for any $s\ge 0$ and $x>0$ we have
 \beqnn
\mbf{P}[f(\xi(r):r\in[0,s])|\tau_{-x}>t]
 \to
\mbf{Q}_x[f(\xi(r)-x: r\in [0,s])],\qquad t\to \infty.
 \eeqnn
\end{proposition}

The key of the proof of Theorem~\ref{th2.9}(2) is the observation that the asymptotics of the expectation (\ref{eq1.2}) only depends on the sample paths of the L\'evy process with slowly decreasing local infimum so that we can use the above two propositions to determine the limiting coefficient. To show clearly the main ideas of the proof, we write the main steps into a series of lemmas.

\begin{lemma}\label{th3.3} There exists a constant $C\ge 0$ such that
 \beqnn
\limsup_{t\to \infty} t^{1/2} \mbf{P}[A_t^\alpha(\xi)^{-\beta/\alpha}]\le C.
 \eeqnn
\end{lemma}

\proof By Theorem~A in Kozlov (1976), there exists a constant $C=C_\beta\ge 0$ such that, as $t\to \infty$,
 \beqnn
\mbf{P}\Big[\exp\Big\{ \min_{k\le [t]-1} \beta \xi(k)\Big\}\Big]
 \sim
C([t]-1)^{-1/2}\sim Ct^{-1/2}.
 \eeqnn
Then the desired result follows from Lemma~\ref{th2.3}. \qed

\begin{lemma}\label{th3.4} Suppose that Condition~\ref{th2.8} holds. Then for any $x>0$ there is a constant $C=C_x\ge 0$ so that
 \beqnn
\mbf{P}[e^{-\alpha \xi(r)};\tau_{-x}> t]
 \le
C(t-r)^{-1/2} r^{-3/2}, \qquad t>r>0.
 \eeqnn
\end{lemma}

\proof By the Markov property and the spatial homogeneity of the L\'evy process,
 \beqnn
\mbf{P}[e^{-\alpha \xi(r)};\tau_{-x}> t]
 \ar=\ar
\mbf{P}[e^{-\alpha \xi(r)}\mbf{P}_{\xi(r)}(\tau_{-x}>t-r);\tau_{-x}> r]\ccr
 \ar=\ar
\mbf{P}[e^{-\alpha \xi(r)}\mbf{P}(\tau_{-x-\xi(r)}>t-r);\tau_{-x}> r].
 \eeqnn
By Corollary~5.3 in Kyprianou (2014, p.118) we have $\hat{V}(y)\sim y/\hat{\mbf{P}}[H(1)]$ as $y\to \infty$. Then the function $y\mapsto e^{-\lambda y}\hat{V}(y)$ is bounded on $[0,\infty)$ for any $\lambda\in (0,\alpha)$. For $x>x_0>0$ we can use Proposition~\ref{th3.1} to see
 \beqnn
\mbf{P}[e^{-\alpha \xi(r)};\tau_{-x}> t]
 \ar=\ar
\mbf{P}[e^{-\alpha \xi(r)}\mbf{P}(\tau_{-x-\xi(r)}> t-r); \tau_{-x}> r]\ccr
 \ar=\ar
\mbf{P}[e^{-\alpha \xi(r)}\mbf{P}(\tau_{-x-\xi(r)}> t-r); \xi(r)\le -x_0,\tau_{-x}> r]\ccr
 \ar\ar
+\,\mbf{P}[e^{-\alpha \xi(r)}\mbf{P}(\tau_{-x-\xi(r)}> t-r); \xi(r)>-x_0,\tau_{-x}> r]\ccr
 \ar\le\ar
\mbf{P}[e^{-\alpha \xi(r)}\mbf{P}(\tau_{x_0-x}> t-r); \tau_{-x}> r]\ccr
 \ar\ar
+\,C(x)(t-r)^{-1/2}\mbf{P}[e^{-\alpha \xi(r)}\hat{V}(x+\xi(r)); \xi(r)>-x_0,\tau_{-x}> r]\ccr
 \ar\le\ar
C(x)\hat{V}(x-x_0)(t-r)^{-1/2}\mbf{P}[e^{-\alpha \xi(r)}; \tau_{-x}> r]\ccr
 \ar\ar
+\, C(x)(t-r)^{-1/2}\mbf{P}[e^{-(\alpha-\lambda)\xi(r)}e^{\lambda x}; \xi(r)>-x_0,\tau_{-x}> r]\ccr
 \ar\le\ar
C(x)\hat{V}(x-x_0)(t-r)^{-1/2}e^{\alpha x}\hat{V}(x) r^{-3/2}\ccr
 \ar\ar
+\, C(x)(t-r)^{-1/2}e^{\lambda x}\mbf{P}[e^{-(\alpha-\lambda) \xi(r)}; \tau_{-x}> r]\ccr
 \ar\le\ar
C(x)[\hat{V}(x-x_0)+1]e^{\alpha x}\hat{V}(x)(t-r)^{-1/2} r^{-3/2}.
 \eeqnn
That gives the desired result. \qed

\begin{lemma}\label{th3.5} Suppose that $F$ is a bounded function satisfying Condition~\ref{th2.6}. Then there is a constant $C\ge 0$ so that, for any $x>0$,
 \beqnn
\limsup_{t\to \infty}t^{1/2}\mbf{P}[F(A_t^\alpha(\xi));\tau_{-x}\le t]
 \le
Ce^{-\beta x}[1+\hat{V}(x)].
 \eeqnn
\end{lemma}

\proof Since $F$ is bounded and satisfies Condition~\ref{th2.6}, there is a constant $C_1\ge 0$ such that $F(z)\le C_1(1\land z^{-\beta/\alpha})$ for all $z>0$. By Lemma~\ref{th3.3} we can find an integer $t_0\ge 3$ and some constant $C\ge 0$ such that
 \beqnn
t^{1/2}\mbf{P}[A_t^\alpha(\xi)^{-\beta/\alpha}]\le C, \qquad t\ge t_0.
 \eeqnn
On the hand, by Proposition~\ref{th3.1}(1), for any $\epsilon>0$ we have, as $t\to \infty$,
 \beqnn
\mbf{P}(\tau_{-x}>t)^{-1}\mbf{P}(t<\tau_{-x}\le t+\epsilon)
 \ar=\ar
\mbf{P}(\tau_{-x}>t)^{-1}[\mbf{P}(\tau_{-x}>t)-\mbf{P}(\tau_{-x}>t+\epsilon)]\ccr
 \ar=\ar
1-\mbf{P}(\tau_{-x}>t)^{-1}\mbf{P}(\tau_{-x}>t+\epsilon) \ccr
 \ar\sim\ar
1-t^{1/2}(t+\epsilon)^{-1/2}\sim 2\epsilon t^{-1}.
 \eeqnn
It follows that
 \beqlb\label{eq3.1}
\mbf{P}(t<\tau_{-x}\le t+\epsilon)\sim\frac{1}{\sqrt{2\pi{\it\Phi}''(0)}}\hat{\mbf{P}}[H(1)] \hat{V}(x-) \epsilon t^{-3/2}.
 \eeqlb
By the strong Markov property, up to some adjustments of value of $C\ge 0$, we have
 \beqnn
\mbf{P}[F(A_t^\alpha(\xi));\tau_{-x}\le t]
 \ar\le\ar
C_1\mbf{P}\Big[1\land \Big(\int_0^t e^{-\alpha \xi(r)}dr\Big)^{-\beta/\alpha};\tau_{-x}\le t\Big] \cr
 \ar\le\ar
C_1\sum_{i=1}^{[t]-t_0}\mbf{P}\Big[\Big(\int_{\tau_{-x}}^{t+\tau_{-x}-i} e^{-\alpha \xi(r)}dr\Big)^{-\beta/\alpha};i-1<\tau_{-x}\le i\Big] \cr
 \ar\ar\qquad\qquad\qquad
+\, C_1\mbf{P}([t]-t_0<\tau_{-x}\le t) \cr
 \ar\le\ar
C\sum_{i=1}^{[t]-t_0}\mbf{P}\Big[e^{\beta \xi(\tau_{-x})}\mbf{P}_{\xi(\tau_{-x})}[(A_{t-i}^\alpha(\xi)^{-\beta/\alpha}]; i-1<\tau_{-x}\le i\Big] \cr
 \ar\ar\qquad\qquad\qquad
+\, C\hat{V}(x)(t_0+1)([t]-t_0)^{-3/2} \cr
 \ar\le\ar
Ce^{-\beta x}\sum_{i=1}^{[t]-t_0}(t-i)^{-1/2}\mbf{P}(i-1<\tau_{-x}\le i) \cr
 \ar\ar\qquad\qquad\qquad
+\, C\hat{V}(x)(t_0+1)([t]-t_0)^{-3/2} \cr
 \ar\le\ar
Ce^{-\beta x}(t-1)^{-1/2} + C\hat{V}(x) e^{-\beta x}\sum_{i=2}^{[t]-t_0}(t-i)^{-1/2}(i-1)^{-3/2} \cr
 \ar\ar\qquad\qquad\qquad
+\, C\hat{V}(x)(t_0+1)([t]-t_0)^{-3/2}.
 \eeqnn
Observe that
 \beqnn
\sum_{i=2}^{[t]-t_0}(t-i)^{-1/2}(i-1)^{-3/2}
 \ar\le\ar
(t-2)^{-1/2} + \int_3^t(t-s)^{-1/2}(s-2)^{-3/2}ds\cr
 \ar\le\ar
(t-2)^{-1/2} + (t/2)^{-1/2} \int_3^{t/2}(s-2)^{-3/2}ds\cr
 \ar\ar
+(t/2-2)^{-3/2}\int_{t/2}^t(t-s)^{-1/2} ds\cr
 \ar\le\ar
(t-2)^{-1/2} + (t/2)^{-1/2} \int_3^{\infty}(s-2)^{-3/2}ds\cr
 \ar\ar
+(t/2-2)^{-3/2}\int_0^t(t-s)^{-1/2} ds\cr
 \ar\le\ar
(t-2)^{-1/2} + 2(t/2)^{-1/2}+2(t/2-2)^{-3/2}t^{1/2}.
 \eeqnn
By combining the above estimates we get desired result. \qed

\begin{lemma}\label{th3.6} Suppose that Condition~\ref{th2.8} holds and $F$ is a globally Lipschitz function on $(0,\infty)$. Then for any $x>0$ we have
 \beqnn
\lim_{s\to \infty}\limsup_{t\to \infty} t^{1/2} \mbf{P}[|F(A_s^\alpha(\xi)) - F(A_t^\alpha(\xi))|; \tau_{-x}> t] = 0.
 \eeqnn
\end{lemma}

\proof Since $F$ is decreasing and globally Lipschitz, there exists a constant $C>0$ such that
 \beqnn
0\le F(A_s^\alpha(\xi)) - F(A_t^\alpha(\xi))\le C \int_s^t e^{-\alpha \xi(r)}dr, \qquad t\ge s\ge t.
 \eeqnn
Then it suffices to prove
 \beqlb\label{eq3.2}
 \lim_{s\to \infty}\lim_{t\to \infty} t^{1/2} \int_s^t \mbf{P}[e^{-\alpha \xi(r)};\tau_{-x}> t]dr =0.
 \eeqlb
By Lemma~\ref{th3.4}, for any $s>1$ we have
 \beqnn
\lefteqn{\limsup_{t\to \infty} t^{1/2} \int_s^t \mbf{P}[e^{-\alpha \xi(r)};\tau_{-x}> t]dr}\qquad\ar\ar\cr
 \ar\le\ar
\limsup_{t\to \infty} Ct^{1/2}\int_s^t(t-r)^{-1/2}r^{-3/2}dr\cr
 \ar\le\ar
\limsup_{t\to \infty} C\Big(\int_s^{t/2}r^{-3/2}dr+ t^{-1}\int_{t/2}^t(t-r)^{-1/2}dr\Big)\cr
 \ar\le\ar
\limsup_{t\to \infty} C\Big(\int_s^{\infty}r^{-3/2}dr+ t^{-1}\int_{0}^t(t-r)^{-1/2}dr\Big)\le Cs^{-1/2}.
 \eeqnn
That proves (\ref{eq3.2}). \qed

\begin{lemma}\label{th3.7} Suppose that Condition~\ref{th2.8} holds and $F$ is a globally Lipschitz function on $(0,\infty)$. Then for any $x>0$ we have
 \beqlb\label{eq3.3}
\lim_{t\to \infty} t^{1/2}\mbf{P}[F(A_t^\alpha(\xi))]);\tau_{-x}> t]
 =
\sqrt{\frac{2}{\pi{\it\Phi}''(0)}}\hat{\mbf{P}}[H(1)]D_2(x,\alpha,F),
 \eeqlb
where
 \beqlb\label{eq3.4}
D_2(x,\alpha,F)
 =
\hat{V}(x-)\mbf{Q}_x[F(e^{-\alpha x}A_\infty^\alpha(\xi))].
 \eeqlb
Furthermore, the function $x\mapsto D_2(x,\alpha,F)$ on $(0,\infty)$ is increasing, strictly positive and bounded.
\end{lemma}

\proof From Lemma~\ref{th2.4} it follows that $\mbf{Q}_x(A_\infty^\alpha(\xi)< \infty) = 1$. Since $F(z)>0$ for each $z>0$, we have $D_2(x,\alpha,F)> 0$. By the spacial homogeneity of the L\'evy process and Proposition~\ref{th3.2}, we have
 \beqnn
\lim_{t\to \infty} \mbf{P}[F(A_s^\alpha(\xi))|\tau_{-x}> t]
 \ar=\ar
\lim_{t\to \infty} \mbf{P}_x[F(A_s^\alpha(\xi-x))|\tau_0> t] \ccr
 \ar=\ar
\lim_{t\to \infty} \mbf{P}_x[F(e^{-\alpha x}A_s^\alpha(\xi))|\tau_0> t]
 =
\mbf{Q}_x[F(e^{-\alpha x}A_s^\alpha(\xi))].
 \eeqnn
Combining this with Proposition~\ref{th3.1}(1) and Lemma~\ref{th3.6},
 \beqnn
\lim_{t\to \infty} t^{1/2}\mbf{P}[F(A_t^\alpha(\xi));\tau_{-x}> t]
 \ar=\ar
\lim_{s\to \infty} \lim_{t\to \infty} t^{1/2}\mbf{P}[F(A_s^\alpha(\xi));\tau_{-x}> t] \ccr
 \ar=\ar
\lim_{s\to \infty} \lim_{t\to \infty} t^{1/2} \mbf{P}(\tau_{-x}> t)
\mbf{P}[F(A_s^\alpha(\xi))|\tau_{-x}> t] \ccr
 \ar=\ar
\lim_{s\to \infty} \sqrt{\frac{2}{\pi{\it\Phi}''(0)}}\hat{\mbf{P}}[H(1)]\hat{V}(x-) \mbf{Q}_x[F(e^{-\alpha x}A_s^\alpha(\xi))] \cr
 \ar=\ar
\sqrt{\frac{2}{\pi{\it\Phi}''(0)}}\hat{\mbf{P}}[H(1)] \hat{V}(x-)\mbf{Q}_x[F(e^{-\alpha x}A_\infty^\alpha(\xi))].
 \eeqnn
Since left-hand side of (\ref{eq3.3}) is increasing in $x>0$, so is $D_2(x,\alpha,F)$. By Lemma~\ref{th3.3} this function is bounded on $(0,\infty)$. \qed

\medskip

\noindent\textit{Proof of Theorem~\ref{th2.9}(2).~} We first consider the special case where $F$ is globally Lipschitz and hence bounded on $(0,\infty)$. By Lemma~\ref{th3.5} we have
 \beqnn
\lim_{x\to \infty}\limsup_{t\to \infty} t^{1/2}
\mbf{P}[F(A_t^\alpha(\xi)); \tau_{-x}\le t] = 0.
 \eeqnn
Then we can use Lemma~\ref{th3.7} to see
 \beqnn
\lim_{t\to \infty}t^{1/2}\mbf{P}[F(A_t^\alpha(\xi))]
 \ar=\ar
\lim_{x\to \infty}\lim_{t\to \infty}t^{1/2}\mbf{P}[F(A_t^\alpha(\xi));\tau_{-x}>t]\cr
 \ar\ar\qquad
+ \lim_{x\to \infty} \lim_{t\to \infty} t^{1/2}\mbf{P}[F(A_t^\alpha(\xi));\tau_{-x}\le t] \cr
 \ar=\ar
\lim_{x\to \infty} \sqrt{\frac{2}{\pi{\it\Phi}''(0)}}\hat{\mbf{P}}[H(1)]D_2(x,\alpha,F).
 \eeqnn
By Lemma~\ref{th3.7}, the limit $D_2(\alpha,F) := \lim_{x\to \infty} D_2(x,\alpha,F)$ is finite, strictly positive and given by (\ref{eq2.7}). Then the result follows in the special case. In the general case, for $n\ge 1$ let $F_n(y) = F(1/n)1_{\{y\le 1/n\}} + F(y)1_{\{y>1/n\}}$ and $G_n(y) = F(y)-F_n(y)$. Then each $F_n$ is globally Lipschitz, so $D_2(x,\alpha,F_n)$ and $D_2(\alpha,F_n)$ can be defined. Clearly, the limit $D_2(x,\alpha,F) := \lim_{n\to \infty} D_2(x,\alpha,F_n)$ exists and is given by (\ref{eq3.4}). It is also easy to see that $D_2(x,\alpha,F)$ is bounded, strictly positive and increasing on $(0,\infty)$. Then the limit $D_2(\alpha,F) := \lim_{x\to \infty} D_2(x,\alpha,F)$ exists and it is finite and strictly positive. Observe that
 \beqnn
\mbf{P}[F(A_t^\alpha(\xi))]
 =
\mbf{P}[G_n(A_t^\alpha(\xi))]
 +
\mbf{P}[F_n(A_t^\alpha(\xi))].
 \eeqnn
By multiplying this by $t^{1/2}$ and taking the limit we have
 \beqlb\label{eq3.5}
\lim_{t\to \infty}t^{1/2}\mbf{P}[F(A_t^\alpha(\xi))]
 =
\lim_{t\to \infty}t^{1/2}\mbf{P}[G_n(A_t^\alpha(\xi))]
+ \sqrt{\frac{2}{\pi{\it\Phi}''(0)}}\hat{\mbf{P}}[H(1)] D_2(\alpha,F_n).
 \eeqlb
Under our general assumption, we can find constants $C\ge 0$ and $\beta_1\in \mcr{D}_+^\circ({\it\Phi})\cap (\beta_0,\infty)$ so that
 \beqnn
\lim_{t\to \infty}t^{1/2}\mbf{P}[G_n(A_t^\alpha(\xi))]
 \ar\le\ar
C\limsup_{t\to \infty}t^{1/2}\mbf{P}[A_t^\alpha(\xi)^{-\beta_0/\alpha}; A_t^\alpha(\xi)\le 1/n] \ccr
 \ar\le\ar
Cn^{-(\beta_1-\beta_0)/\alpha}\limsup_{t\to \infty}t^{1/2}\mbf{P}[A_t^\alpha(\xi)^{-\beta_1/\alpha}] \ccr
 \ar\le\ar
Cn^{-(\beta_1-\beta_0)/\alpha}.
 \eeqnn
Then for $k\ge n\ge 1$ we have
 \beqnn
0\ar\le\ar D_2(x,\alpha,F_k)-D_2(x,\alpha,F_n) \ccr
 \ar=\ar
\sqrt{\frac{\pi{\it\Phi}''(0)}{2}}\hat{\mbf{P}}[H(1)]^{-1} \lim_{t\to \infty}t^{1/2}\mbf{P}[F_k(A_t^\alpha(\xi))-F_n(A_t^\alpha(\xi)), \tau_{-x}>t] \cr
 \ar\le\ar
\sqrt{\frac{\pi{\it\Phi}''(0)}{2}}\hat{\mbf{P}}[H(1)]^{-1} \lim_{t\to \infty}t^{1/2}\mbf{P}[G_n(A_t^\alpha(\xi)), \tau_{-x}>t] \cr
 \ar\le\ar
\sqrt{\frac{\pi{\it\Phi}''(0)}{2}}\hat{\mbf{P}}[H(1)]^{-1} Cn^{-(\beta_1-\beta_0)/\alpha}.
 \eeqnn
By letting $k\to \infty$ in the above we see
 \beqnn
0\le D_2(x,\alpha,F)-D_2(x,\alpha,F_n)
 \le
\sqrt{\frac{\pi{\it\Phi}''(0)}{2}}\hat{\mbf{P}}[H(1)]^{-1} Cn^{-(\beta_1-\beta_0)/\alpha},
 \eeqnn
and hence
 \beqnn
0\le D_2(\alpha,F)-D_2(\alpha,F_n)
 \le
\sqrt{\frac{\pi{\it\Phi}''(0)}{2}}\hat{\mbf{P}}[H(1)]^{-1} Cn^{-(\beta_1-\beta_0)/\alpha}.
 \eeqnn
Then we can let $n\to \infty$ in (\ref{eq3.5}) to get the result. \qed

We remark that for the characteristic exponent given by (\ref{eq2.12}) a result similar to Proposition~\ref{th3.1}(2) was established in Lemma~A-(b) of Hirano (2001). One may check that all the arguments given above carry over to that case under obvious modifications.

\section{Transient L\'evy processes}\label{StronglySubcritical}

\setcounter{equation}{0}

In this section, we give the proof of Theorem~\ref{th2.9} in the case where $\xi$ is a transient L\'evy process, i.e., $\mbf{P}[\xi(1)]\neq 0$. In fact, the result of regime (1) is a simple consequence of Proposition~\ref{th2.1}. The proof for regime (3) is based on an extension of Theorem~2 of Hirano (2001), where the exponential functional was approximated by functionals of two independent processes obtained by transformations. The proofs for regimes (4) and (5) are also based on transformations of the underlying L\'evy process.

\medskip

\noindent\textit{Proof of Theorem~\ref{th2.9}(1).~} Since $\mbf{P}[\xi(1)]> 0$, we have $\lim_{t\to \infty}\xi(t) = \infty$ and hence $\mbf{P}(A_\infty^\alpha(\xi)< \infty) = 1$ by Proposition~\ref{th2.1}. Then the result is immediate for any bounded function $F$. For an unbounded function $F$, under the general assumption we have
 \beqnn
\mbf{P}[F(A_t^\alpha(\xi))]
 \le
C_0\mbf{P}[A_t^\alpha(\xi)^{-\beta_0/\alpha}] + \mbf{P}[F(A_t^\alpha(\xi))1_{\{A_t^\alpha(\xi)\ge 1\}}].
 \eeqnn
The right-hand side is finite by Lemma~\ref{th2.2}. By monotone convergence, we get $\mbf{P}[F(A_t^\alpha(\xi))]\to \mbf{P}[F(A_\infty^\alpha(\xi))]$ decreasingly as $t\to \infty$. Clearly, the limit is finite and strictly positive. \qed

\begin{proposition}\label{th4.1} Suppose that Condition~\ref{th2.8} holds and there exists $\varrho\in \mcr{D}_+^\circ({\it\Phi})$ satisfying ${\it\Phi}'(\varrho) = 0$. Then for any $x>0$ and $\theta\in \mcr{D}({\it\Phi})\cap (-\varrho,\infty)$ we have, as $t\to \infty$,
 \beqnn
\mbf{P}(\tau_{-x}>t)
 \sim
Me^{\varrho x}\hat{V}^{(\varrho)}(x-) t^{-3/2} e^{{\it\Phi}(\varrho)t} \int_{0}^\infty e^{-\varrho z}V^{(\varrho)}(z)dz
 \eeqnn
and
 \beqnn
\mbf{P}(e^{-\theta \xi(t)};\tau_{-x}>t)
 \sim
M\hat{V}^{(\varrho)}(x-)e^{(\theta+\varrho)x}t^{-3/2}e^{{\it\Phi}(\varrho)t} \int_{0}^\infty e^{-(\theta+\varrho)z}V^{(\varrho)}(z)dz,
 \eeqnn
where
 \beqnn
M = \frac{1}{\sqrt{2\pi{\it\Phi}''(\varrho)}}\exp\Big\{\int_0^\infty (e^{-t}-1)t^{-1}e^{-t{\it\Phi}(\varrho)} \mbf{P}(\xi(t)=0)dt\Big\}.
 \eeqnn
\end{proposition}

\proof We only need to show the second result since the first one is its special case with $\theta = 0$. By the definition of $\mbf{P}^{(\varrho)}$ we have
 \beqnn
e^{-{\it\Phi}(\varrho)t}\mbf{P}(e^{-\theta \xi(t)};\tau_{-x}>t)
 =
\mbf{P}^{(\varrho)}(e^{-(\theta+\varrho) \xi(t)};\tau_{-x}>t).
 \eeqnn
Since $\mbf{P}^{(\varrho)}[\xi(1)] = {\it\Phi}'(\varrho) = 0$ and $\mbf{P}^{(\varrho)}[\xi(1)^2] = {\it\Phi}''(\varrho)< \infty$, the desired result follows by Proposition~\ref{th3.1}(2). \qed

Recall that $D[0,s]$ denotes the space of c\`adl\`ag real functions on $[0,s]$ equipped with the Skorokhod topology. Let $C_0(\mbb{R}_+)$ and $C_0(\mbb{R}_+^2)$ denote respectively the spaces of continuous function on $\mbb{R}_+$ and $\mbb{R}_+^2$ vanishing at infinity. The following proposition is a simple extension of Theorem~2-(a) in Hirano (2001).

\begin{proposition}\label{th4.2} Suppose that Condition~\ref{th2.8} holds and there exists $\varrho\in \mcr{D}_+^\circ({\it\Phi})$ satisfying ${\it\Phi}'(\varrho) = 0$. Let $H\in C_0(\mbb{R}_+^2)$ and let $f,g$ be continuous functions on $D[0,s]$. Then for any $x>0$ we have
 \beqlb\label{eq4.1}
\lefteqn{\lim_{t\to \infty} \mbf{P}\Big[H(f((\xi(r))_{0\le r\le s}),g((\xi(t-r))_{0\le r\le s})\big|\tau_{-x}>t\Big]}\qquad\ar\ar\cr
 \ar=\ar
\frac{1}{h(\varrho)}\int_0^\infty e^{-\varrho y}V^{(\varrho)}(y) \mbf{Q}_{(x,y)}^{(\varrho)} [H(f((\xi(r)-x)_{0\le r\le s}), g((\hat{\xi}(r)-x)_{0\le r\le s})]dy,
 \eeqlb
where
 \beqnn
h(\varrho) = \int_0^\infty e^{-\varrho y}V^{(\varrho)}(y) dy.
 \eeqnn
\end{proposition}

\proof If $H(x_1,x_2) = G_1(x_1)G_2(x_2)$ for $G_1,G_2\in C_0(\mbb{R}_+)$, we have (\ref{eq4.1}) by Theorem~2-(a) in Hirano (2001). The general result follows by the Stone-Weierstrass Theorem. \qed

\begin{lemma}\label{th4.3} Suppose that Condition~\ref{th2.8} holds and there exists $\varrho\in \mcr{D}_+^\circ({\it\Phi})$ satisfying ${\it\Phi}'(\varrho) = 0$. Then for any $\beta\in \mcr{D}({\it\Phi})\cap (\varrho,\infty)$ there exists a constant $C=C_\beta> 0$ such that
 \beqnn
\limsup_{t\to \infty} t^{3/2}e^{-{\it\Phi}(\varrho)t} \mbf{P}[A_t^\alpha(\xi)^{-\beta/\alpha}]\le C.
 \eeqnn
\end{lemma}

\proof By Lemma~7(3) in Hirano (1998), under the conditions of the lemma there exists a constant $\eta(\beta)\ge 0$ such that
 \beqnn
\mbf{P}\Big[\exp\Big\{ \min_{k\le [t]-1}\beta \xi(k)\Big\}\Big]
 \le
\eta(\beta)([t]-1)^{-3/2}e^{{\it\Phi}(\varrho)([t]-1)}
 \le
2^{3/2}\eta(\beta)e^{-{\it\Phi}(\varrho)}t^{-3/2}e^{{\it\Phi}(\varrho)t}.
 \eeqnn
Then we have the result by Lemmas~\ref{th2.2} and~\ref{th2.3}. \qed

\begin{lemma}\label{th4.4} Suppose that Conditions~\ref{th2.5} and~\ref{th2.8} hold and there exists $\varrho\in \mcr{D}_+^\circ({\it\Phi})$ satisfying ${\it\Phi}'(\varrho) = 0$. If $F$ is a bounded function, in addition, then for any $x>0$ and $\beta\in \mcr{D}({\it\Phi})\cap (\varrho,\infty)$ there is a constant $C\ge 0$ so that
 \beqnn
\limsup_{t\to \infty} t^{3/2}e^{-t{\it\Phi}(\varrho)} \mbf{P}[F(A_t^\alpha(\xi));\tau_{-x}\le t-\epsilon]
 \le
Ce^{-\beta x} + Ce^{-(\beta-\varrho)x}\hat{V}^{(\varrho)}(x),
 \eeqnn
and hence
 \beqnn
\lim_{x\to \infty}\limsup_{t\to \infty} t^{3/2}e^{-t{\it\Phi}(\varrho)} \mbf{P}[F(A_t^\alpha(\xi)); \tau_{-x}\le t-\epsilon] = 0.
 \eeqnn
\end{lemma}

\proof This is a modification of the proof of Lemma~\ref{th3.5}. By Lemma~\ref{th4.3}, there exists a constant $C\ge 0$ so that
 $$
\mbf{P}[A_t^\alpha(\xi)^{-\beta/\alpha}]\le Ct^{-3/2}e^{{\it\Phi}(\varrho)t}, \qquad t\ge 2.
 $$
For any $\eta>0$, one can use Proposition~\ref{th4.1} to see as the derivation of (\ref{eq3.1}) that, as $t\to \infty$,
 \beqnn
\mbf{P}(t<\tau_{-x}\le t+\eta)
 \sim
Me^{\varrho x}\hat{V}^{(\varrho)}(x-) (1-e^{{\it\Phi}(\varrho)\eta}) t^{-3/2} e^{{\it\Phi}(\varrho)t} \int_{0}^\infty e^{-\varrho z}V^{(\varrho)}(z)dz,
 \eeqnn
where
 \beqnn
M = \frac{1}{\sqrt{2\pi{\it\Phi}''(\varrho)}}\exp\Big\{\int_0^\infty (e^{-t}-1)t^{-1}e^{-t{\it\Phi}(\varrho)} \mbf{P}(\xi(t)=0)dt\Big\}.
 \eeqnn
By adjusting the value of $C\ge 0$ we have, for $t\ge 3$ and $0<\epsilon \le 1$,
 \beqnn
\mbf{P}[F(A_t^\alpha(\xi));\tau_{-x}\le t-\epsilon]
 \ar\le\ar
C \sum_{i=1}^{[t]-2}\mbf{P}\Big[\Big(\int_{\tau_{-x}}^{t-i+\tau_{-x}} e^{-\alpha \xi(r)}dr\Big)^{-\beta/\alpha};i-1<\tau_{-x}\le i\Big] \cr
 \ar\ar
+\, C\mbf{P}\Big[\Big(\int_{\tau_{-x}}^{\tau_{-x}+\epsilon} e^{-\alpha \xi(r)}dr\Big)^{-\beta/\alpha}; [t]-2<\tau_{-x}\le t-\epsilon\Big] \ccr
 \ar=\ar
C\sum_{i=1}^{[t]-2}\mbf{P}\big\{e^{-\beta \xi(\tau_{-x})} \mbf{P}_{\xi(\tau_{-x})}[A_{t-i}^\alpha(\xi)^{-\beta/\alpha}]; i-1<\tau_{-x}\le i\big\}\ccr
 \ar\ar
+\, C\mbf{P}\big\{e^{-\beta \xi(\tau_{-x})}\mbf{P}_{\xi(\tau_{-x})}[A_\epsilon^\alpha(\xi)^{-\beta/\alpha}]; [t]-2<\tau_{-x}\le t\big\} \ccr
 \ar\le\ar
Ce^{-\beta x}\sum_{i=1}^{[t]-2}\mbf{P}\big\{\mbf{P}_{\xi(\tau_{-x})}[A_{t-i}^\alpha(\xi)^{-\beta/\alpha}]; i-1<\tau_{-x}\le i\big\}\ccr
 \ar\ar
+\, Ce^{-\beta x}\mbf{P}\big\{\mbf{P}_{\xi(\tau_{-x})}[A_\epsilon^\alpha(\xi)^{-\beta/\alpha}]; [t]-2<\tau_{-x}\le t\big\} \ccr
 \ar\le\ar
Ce^{-\beta x}\sum_{i=1}^{[t]-2}e^{{\it\Phi}(\varrho)(t-i)}(t-i)^{-3/2} \mbf{P}(i-1<\tau_{-x}\le i) \ccr
 \ar\ar\quad
+\, Ce^{-\beta x}\mbf{P}([t]-2<\tau_{-x}\le t) \ccr
 \ar\le\ar
C e^{-(\beta-\varrho)x} \hat{V}^{(\varrho)}(-x) e^{{\it\Phi}(\varrho)(t-1)} \sum_{i=2}^{[t]-2}(t-i)^{-3/2} (i-1)^{-3/2} \ccr
 \ar\ar
+\, Ce^{-\beta x}e^{{\it\Phi}(\varrho)(t-1)}(t-1)^{-3/2} \ccr
 \ar\ar
+\, Ce^{-(\beta-\varrho)x}\hat{V}^{(\varrho)}(-x) (1-e^{3{\it\Phi}(\varrho)}) ([t]-2)^{-3/2}e^{{\it\Phi}(\varrho)([t]-2)},
 \eeqnn
where
 \beqnn
\sum_{i=2}^{[t]-2}(t-i)^{-3/2}(i-1)^{-3/2}
 \ar\le\ar
(t-2)^{-3/2} + \int_2^{t-2}(t-s-1)^{-3/2}(s-1)^{-3/2} ds\cr
 \ar\le\ar
(t-2)^{-3/2} + \Big(\frac{t}{2}-1\Big)^{-3/2} \int_2^{t/2}(s-1)^{-3/2}ds\cr
 \ar\ar
+ \Big(\frac{t}{2}-1\Big)^{-3/2}\int_{t/2}^{t-2}(t-s-1)^{-3/2} ds\cr
 \ar\le\ar
(t-2)^{-3/2} + 4\Big(\frac{t}{2}-1\Big)^{-3/2}.
 \eeqnn
Then we get the desired result. \qed

\begin{lemma}\label{th4.5} Suppose that Conditions~\ref{th2.5} and~\ref{th2.8} hold and there exists $\varrho\in \mcr{D}_+^\circ({\it\Phi})$ satisfying ${\it\Phi}'(\varrho) = 0$. If $F$ is a bounded function, in addition, then for any $x>0$ we have
 $$
\lim_{s\to \infty}\limsup_{t\to \infty} t^{3/2}e^{-t{\it\Phi}(\varrho)} \mbf{P}\bigg[F\Big(\int_{[0,s]\cup [t-s,t]} e^{-\alpha \xi(r)}dr\Big) -F(A_t^\alpha(\xi));\tau_{-x}> t\bigg]=0.
 $$
\end{lemma}

\proof As in the proof of Lemma~\ref{th3.4}, one can use Proposition~\ref{th4.1} to see there is a constant $C=C_x\ge 0$ so that
 \beqnn
e^{-{\it\Phi}(\varrho)t}\mbf{P}\Big[e^{-\alpha \xi(r)};\tau_{-x}> t\Big]\le C(t-r)^{-3/2} r^{-3/2}.
 \eeqnn
By elementary analysis we have
 \beqnn
\int_s^{t-s}(t-r)^{-3/2} r^{-3/2}dr\le 2 (t/2)^{-3/2}\int_s^\infty r^{-3/2}dr= 8\sqrt{2}t^{-3/2}s^{-1/2}.
 \eeqnn
Then we can prove the statement as in the proof of Lemma~\ref{th3.6}. \qed

\begin{lemma}\label{th4.6} Suppose that Conditions~\ref{th2.5} and~\ref{th2.8} hold and there exists $\varrho\in \mcr{D}_+^\circ({\it\Phi})$ satisfying ${\it\Phi}'(\varrho) = 0$. If $F$ is a bounded function, in addition, then for any $x>0$ we have
 \beqlb\label{eq4.2}
\lim_{t\to\infty} t^{3/2}e^{-t{\it\Phi}(\varrho)}\mbf{P}[F(A_t^\alpha(\xi));\tau_{-x}> t]
 =
\frac{c(\varrho)}{\sqrt{2\pi{\it\Phi}''(\varrho)}} D_3(x,\alpha,F),
 \eeqlb
where
 \beqnn
D_3(x,\alpha,F) = e^{\varrho x}\hat{V}^{(\varrho)}(x-)\int_0^\infty e^{-\varrho y}V^{(\varrho)}(y) G(x,y)dy,
 \eeqnn
where $G(\cdot,\cdot)$ is defined by (\ref{eq2.10}). Moreover, the function $x\mapsto D_3(x,\alpha,F)$ is bounded, increasing and strictly positive in $(0,\infty)$. \end{lemma}

\proof From Lemma~\ref{th2.4} it follows that $\mbf{Q}_{(x,z)}^{(\varrho)} \{A_\infty^\alpha(\xi) + A_\infty^\alpha(\hat{\xi})< \infty\} = 1$. Then $D_3(x,\alpha,F)$ is well-defined and strictly positive. By
Lemma~\ref{th4.5} we have
 \beqnn
\lefteqn{\lim_{t\to \infty} t^{3/2}e^{-t{\it\Phi}(\varrho)}
\mbf{P}[F(A_t^\alpha(\xi));\tau_{-x}> t]}\qquad\ar\ar\ccr
 \ar=\ar
\lim_{s\to\infty} \lim_{t\to \infty} t^{3/2}e^{-t{\it\Phi}(\varrho)}
\mbf{P}\Big[F\Big(\int_{[0,s]\cup [t-s,t]} e^{-\alpha
\xi(r)}dr\Big);\tau_{-x}> t\Big].
 \eeqnn
By applying Propositions~\ref{th4.1} and \ref{th4.2} with $H(x_1,x_2) = F(x_1+x_2)$, we have, for any $s>0$,
 \beqnn
\lefteqn{\lim_{t\to \infty}t^{3/2}e^{-{\it\Phi}(\varrho)t} \mbf{P}\Big[F\Big(\int_0^s e^{-\alpha \xi(r)}dr + \int_{t-s}^t e^{-\alpha \xi(r)}dr\Big); \tau_{-x}> t\Big]}\qquad\ar\ar\cr
 \ar=\ar
M(x)\int_0^\infty e^{-\varrho y}V^{(\varrho)}(y) \mbf{Q}_{(x,y)}^{(\varrho)} \{F(e^{-\alpha x}[A_s^\alpha(\xi) + A_s^\alpha(\hat{\xi})])\}dy,
 \eeqnn
where
 \beqnn
M(x) = \frac{c(\varrho)}{\sqrt{2\pi{\it\Phi}''(\varrho)}} e^{\varrho x}\hat{V}^{(\varrho)}(x-).
 \eeqnn
Then we get (\ref{eq4.2}) by letting $s\to \infty$. As in the proof of Lemma~\ref{th3.7} one can see that the function $x\mapsto D_3(x,\alpha,F)$ is bounded and increasing in $(0,\infty)$. \qed

\medskip

\noindent\textit{Proof of Theorem~\ref{th2.9}(3).~} We here only consider a bounded function $F$. The general case can be treated similarly as in the proof of Theorem~\ref{th2.9}(2). By Lemma~\ref{th4.6} we have
 \beqnn
\liminf_{t\to \infty} t^{3/2}e^{-t{\it\Phi}(\varrho)} \mbf{P}[F(A_t^\alpha(\xi))]
 \ar\ge\ar
\lim_{x\to \infty}\lim_{t\to \infty}t^{3/2}e^{-t{\it\Phi}(\varrho)} \mbf{P}[F(A_t^\alpha(\xi));\tau_{-x}> t] \ccr
 \ar=\ar
\lim_{x\to \infty}\frac{c(\varrho)}{\sqrt{2\pi{\it\Phi}''(\varrho)}} D_3(x,\alpha,F) \cr
 \ar=\ar
\frac{c(\varrho)}{\sqrt{2\pi{\it\Phi}''(\varrho)}} D_3(\alpha,F).
 \eeqnn
On the other hand, for any $\epsilon>0$ and $x<0$ we can write
 \beqnn
t^{3/2}e^{-t{\it\Phi}(\varrho)} \mbf{P}[F(A_t^\alpha(\xi))]
 \ar=\ar
t^{3/2}e^{-t{\it\Phi}(\varrho)} \mbf{P}[F(A_t^\alpha(\xi));\tau_{-x}\le t-\epsilon]\ccr
 \ar\ar\quad
+\, t^{3/2}e^{-t{\it\Phi}(\varrho)} \mbf{P}[F(A_t^\alpha(\xi));\tau_{-x}> t-\epsilon].
 \eeqnn
By Lemma~\ref{th4.4}, the first term on the right-hand side tends to zero as $t\to \infty$ and $x\to -\infty$. Since $z\mapsto F(z)$ is decreasing, we can use Lemma~\ref{th4.6} to see
 \beqnn
\limsup_{t\to \infty} t^{3/2}e^{-t{\it\Phi}(\varrho)} \mbf{P}[F(A_t^\alpha(\xi))]
 \ar\le\ar
\lim_{x\to \infty}\lim_{t\to \infty}t^{3/2}e^{-t{\it\Phi}(\varrho)} \mbf{P}[F(A_{t-\epsilon}^\alpha(\xi));\tau_{-x}> t-\epsilon] \ccr
 \ar=\ar
\lim_{x\to \infty}e^{-\epsilon{\it\Phi}(\varrho)} \frac{c(\varrho)}{\sqrt{2\pi{\it\Phi}''(\varrho)}} D_3(x,\alpha,F) \cr
 \ar=\ar
e^{-\epsilon{\it\Phi}(\varrho)}\frac{c(\varrho)}{\sqrt{2\pi{\it\Phi}''(\varrho)}} D_3(\alpha,F).
 \eeqnn
Since $\epsilon>0$ was arbitrary, we get
 \beqnn
\limsup_{t\to \infty} t^{3/2}e^{-t{\it\Phi}(\varrho)} \mbf{P}[F(A_t^\alpha(\xi))]
 \le
\frac{c(\varrho)}{\sqrt{2\pi{\it\Phi}''(\varrho)}}D_3(\alpha,F).
 \eeqnn
That gives the desired result. \qed

\begin{lemma}\label{th4.7} Suppose that ${\it\Phi}'(\beta) = 0$. Then for any $\theta\in \mcr{D}_+^\circ({\it\Phi})\cap (\beta,\infty)$ we have
 \beqnn
\lim_{t\to\infty} t^{1/2}e^{-t{\it\Phi}(\beta)}\mbf{P}[A_t^\alpha(\xi)^{-\theta/\alpha}] = 0.
 \eeqnn
\end{lemma}

\proof By the definition of $\mbf{P}^{(\beta)}$ and the property of independent increments of $\{\xi(t): t\ge 0\}$ under this probability measure,
 \beqnn
t^{1/2}e^{-{\it\Phi}(\beta)t}\mbf{P}[A_t^\alpha(\xi)^{-\theta/\alpha}]
 \ar=\ar
t^{1/2}e^{-{\it\Phi}(\beta)t} \mbf{P}[A_t^\alpha(\xi)^{-\beta/\alpha} A_t^\alpha(\xi)^{-(\theta-\beta)/\alpha}]\ccr
 \ar=\ar
t^{1/2}\mbf{P}^{(\beta)}\Big[\Big(\int_0^t e^{\beta [\xi(t)-\xi(s)]} ds\Big)^{-\beta/\alpha} A_t^\alpha(\xi)^{-(\theta-\beta)/\alpha}\Big]\cr
 \ar\le\ar
t^{1/2}\mbf{P}^{(\beta)}\Big[\Big(\int_{t/2}^t e^{\beta [\xi(t)-\xi(s)]} ds\Big)^{-\beta/\alpha} A_{t/2}^\alpha(\xi)^{-(\theta-\beta)/\alpha}\Big]\cr
 \ar=\ar
t^{1/2}\mbf{P}^{(\beta)}\Big[\Big(\int_{t/2}^t e^{\beta [\xi(t)-\xi(s)]} ds\Big)^{-\beta/\alpha}\Big] \mbf{P}^{(\beta)} [A_{t/2}^\alpha(\xi)^{-(\theta-\beta)/\alpha}]\ccr
 \ar=\ar
t^{1/2}\mbf{P}^{(\beta)}[A_{t/2}^\alpha(-\xi)^{-\beta/\alpha}] \mbf{P}^{(\beta)} [A_{t/2}^\alpha(\xi)^{-(\theta-\beta)/\alpha}],
 \eeqnn
where we have used the duality relation in the last equality; see, e.g., Lemma~3.4 in Kyprianou (2014, p.77). Since $\{\xi(t): t\ge 0\}$ is a recurrent L\'evy process under $\mbf{P}^{(\beta)}$, the right-hand side tends to zero as $t\to \infty$ by Lemma~\ref{th3.3}. \qed

\medskip

\noindent\textit{Proof of Theorem~\ref{th2.9}(4).~} By the general assumption, for each $y>0$ there is a constant $C=C_y\ge 0$ so that $F(z)\le Cz^{-\beta_0/\alpha}$ for $0<z\le y$. Upon an adjustment of the value of the constants, we may assume $\beta_0\in \mcr{D}_+^\circ({\it\Phi})\cap (\varrho,\infty)$. By Lemma~\ref{th4.7} we have
 \beqnn
\lefteqn{\limsup_{t\to \infty} t^{1/2}e^{-{\it\Phi}(\beta)t}\mbf{P}[F(A_t^\alpha(\xi)); A_t^\alpha(\xi)< y]}\qquad\ar\ar\ccr
 \ar\le\ar
C\limsup_{t\to \infty} t^{1/2}e^{-{\it\Phi}(\beta)t}\mbf{P}[A_t^\alpha(\xi)^{-\beta_0/\alpha}; A_t^\alpha(\xi)< y]\ccr
 \ar\le\ar
C\limsup_{t\to \infty} t^{1/2}e^{-{\it\Phi}(\beta)t}\mbf{P}[A_t^\alpha(\xi)^{-\beta_0/\alpha}] = 0,
 \eeqnn
and hence
 \beqlb\label{eq4.3}
\lim_{t\to \infty} t^{1/2}e^{-{\it\Phi}(\beta)t} \mbf{P}[F(A_t^\alpha(\xi))]
 =
\lim_{t\to \infty} t^{1/2}e^{-{\it\Phi}(\beta)t} \mbf{P}[F(A_t^\alpha(\xi)); A_t^\alpha(\xi)\ge y].
 \eeqlb
By the duality relation of the L\'evy process,
 \beqnn
A_t^\alpha(\xi)^{-\beta/\alpha}
 \overset{\rm d}=
\Big(\int_0^t e^{\alpha[\xi(t-s)-\xi(t)]}ds\Big)^{-\beta/\alpha}
 \overset{\rm d}=
e^{\beta \xi(t)}A_t^\alpha(-\xi)^{-\beta/\alpha}.
 \eeqnn
It follows that
 \beqlb\label{eq4.4}
e^{-{\it\Phi}(\beta)t}\mbf{P}[A_t^\alpha(\xi)^{-\beta/\alpha}]
 =
e^{-{\it\Phi}(\beta)t}\mbf{P}[e^{\beta \xi(t)}A_t^\alpha(-\xi)^{-\beta/\alpha}]
 =
\mbf{P}^{(\beta)}[A_t^\alpha(-\xi)^{-\beta/\alpha}].
 \eeqlb
By Condition~\ref{th2.7}, given any $\delta\in (0,1)$ we can choose sufficiently large $y>0$ so that
 $$
(1-\delta) Kz^{-\beta/\alpha}\le F(z)\le (1+\delta) Kz^{-\beta/\alpha}, \qquad z>y.
 $$
Then, in view of (\ref{eq4.3}) and (\ref{eq4.4}),
 \beqlb\label{eq4.5}
\lefteqn{\limsup_{t\to \infty} t^{1/2}e^{-{\it\Phi}(\beta)t} \mbf{P}[F(A_t^\alpha(\xi))]}\qquad\ar\ar\ccr
 \ar\le\ar
(1+\delta)K\limsup_{t\to \infty} t^{1/2}e^{-{\it\Phi}(\beta)t} \mbf{P}[A_t^\alpha(\xi)^{-\beta/\alpha}; A_t^\alpha(\xi)\ge y]\ccr
 \ar\le\ar
(1+\delta)K\limsup_{t\to \infty} t^{1/2}e^{-{\it\Phi}(\beta)t} \mbf{P}[A_t^\alpha(\xi)^{-\beta/\alpha}]\ccr
 \ar=\ar
(1+\delta)K\limsup_{t\to \infty} t^{1/2}\mbf{P}^{(\beta)}[A_t^\alpha(-\xi)^{-\beta/\alpha}]\ccr
 \ar=\ar
(1+\delta)K\sqrt{\frac{2}{\pi{\it\Phi}''(\beta)}} \mbf{P}^{(\beta)}[H(1)] D_4(\alpha,\beta),
 \eeqlb
where the last equality follows by Theorem~\ref{th2.9}(2). On the other hand, by the definition of $\mbf{P}^{(\beta)}$ we have
 \beqnn
e^{-{\it\Phi}(\beta)t} \mbf{P}[F(A_t^\alpha(\xi))]
 \ar\ge\ar
(1-\delta)K e^{-{\it\Phi}(\beta)t} \mbf{P}[A_t^\alpha(\xi)^{-\beta/\alpha}; A_t^\alpha(\xi)\ge y]\ccr
 \ar\ge\ar
(1-\delta)K \mbf{P}^{(\beta)}\Big[\Big(\int_0^t e^{\alpha[\xi(t)-\xi(s)]}ds\Big)^{-\beta/\alpha}; A_t^\alpha(\xi)\ge y\Big].
 \eeqnn
Since $\mbf{P}^{(\beta)}[\xi(1)] = {\it\Phi}'(\beta)=0$, from Proposition~\ref{th2.1} it follows that $\mbf{P}^{(\beta)}(A_\infty^\alpha(\xi) = \infty) = 1$. Then, by dominated convergence and the duality relation,
 \beqlb\label{eq4.6}
\lefteqn{\liminf_{t\to \infty} t^{1/2}e^{-{\it\Phi}(\beta)t} \mbf{P}[F(A_t^\alpha(\xi))]}\qquad\ar\ar\ccr
 \ar\ge\ar
(1-\delta)K\liminf_{t\to \infty} t^{1/2}\mbf{P}^{(\beta)}\Big[\Big(\int_0^t e^{\alpha[\xi(t)-\xi(s)]}ds\Big)^{-\beta/\alpha}\Big]\ccr
 \ar=\ar
(1-\delta)K\liminf_{t\to \infty} t^{1/2}\mbf{P}^{(\beta)}[A_t^\alpha(-\xi)^{-\beta/\alpha}]\ccr
 \ar=\ar
(1-\delta)K\sqrt{\frac{2}{\pi{\it\Phi}''(\beta)}} \mbf{P}^{(\beta)}[H(1)] D_4(\alpha,\beta),
 \eeqlb
where we used Theorem~\ref{th2.9}(2) again for the last equality. Since $\delta\in (0,1)$ was arbitrary, we get the desired result by combining (\ref{eq4.5}) and (\ref{eq4.6}). \qed

\medskip

\noindent\textit{Proof of Theorem~\ref{th2.9}(5).~} Under the assumption, we can take $\theta\in \mcr{D}_+^\circ({\it\Phi})\cap (\beta,\infty)$ satisfying ${\it\Phi}(\theta)< {\it\Phi}(\beta)$ and ${\it\Phi}'(\beta)< {\it\Phi}'(\theta)<0$. Then we have
 \beqnn
\limsup_{t\to \infty}e^{-{\it\Phi}(\beta)t}\mbf{P}[A_t^\alpha(\xi)^{-\theta/\alpha}]
 \ar=\ar
\limsup_{t\to \infty}e^{[{\it\Phi}(\theta)-{\it\Phi}(\beta)]t} \mbf{P}^{(\theta)}[A_t^\alpha(-\xi)^{-\theta/\alpha}] \ccr
 \ar=\ar
\limsup_{t\to \infty}e^{[{\it\Phi}(\theta)-{\it\Phi}(\beta)]t} \mbf{P}^{(\theta)}[A_1^\alpha(-\xi)^{-\theta/\alpha}] = 0.
 \eeqnn
The remaining arguments are modifications of those in the proof of Theorem~\ref{th2.9}(4). \qed

\section{Survival probability of the CBRE-process}\label{application}

\setcounter{equation}{0}

Suppose that $({\it\Omega}, \mcr{F},\mcr{F}_t,\mbf{P})$ is a filtered probability space satisfying the usual hypotheses. Let $\sigma\ge 0$ and $b$ be real constants. Let $(z\wedge z^2)\nu(dz)$ be a finite measure on $\mbb{R}$ supported by $\mbb{R}\setminus \{0\}$. Let $\{B(t): t\ge 0\}$ be an $(\mcr{F}_t)$-Brownian motion and $N(ds,dz)$ an $(\mcr{F}_t)$-Poisson random measure on $(0,\infty)\times \mbb{R}$ with intensity $ds\nu(dz)$. Let $\{L(t): t\ge 0\}$ be an $(\mcr{F}_t)$-L\'evy process with the following L\'evy-It\^o decomposition:
 \beqnn
L(t) = \beta t + \sigma B((t) + \int_0^t \int_{[-1,1]}(e^z-1)\tilde{N}(ds,dz) + \int_0^t \int_{[-1,1]^c} (e^z-1) N(ds,dz),
 \eeqnn
where $[-1,1]^c = \mbb{R}\setminus [-1,1]$ and $\tilde{N}(ds,dz) = N(ds,dz) - ds\nu(dz)$. Then $\{L(t): t\ge 0\}$ has with no jump less than $-1$. Let $0<\alpha\le 1$ and let $\{Z_\alpha(t): t\ge 0\}$ be a spectrally positive $(\mcr{F}_t)$-stable process with index $(1+\alpha)$. When $\alpha=1$, we think of $\{Z_\alpha(t): t\ge 0\}$ as a Brownian monition. When $0< \alpha< 1$, we assume $\{Z_\alpha(t): t\ge 0\}$ has L\'{e}vy measure:
 \beqnn
m(dz) = \frac{\alpha 1_{\{z>0\}}dz}{\Gamma(1-\alpha)z^{2+\alpha}}.
 \eeqnn
Suppose that $\{Z_\alpha(t): t\ge 0\}$ and $\{L(t): t\ge 0\}$ are independent of each other. Let $c\ge 0$ be another constant. Let $\{X(t): t\ge 0\}$ be the CBRE-process defined by (\ref{eq1.3}). We can define another L\'evy process $\{\xi(t): t\ge 0\}$ by
 \beqnn
\xi(t) = a_0t + \sigma B(t) + \int_0^t \int_{\mbb{R}} z \tilde{N}(ds,dz),
 \eeqnn
where
 \beqnn
a_0 = \beta-\frac{\sigma^2}{2} - \int_{[-1,1]}(e^z-1-z) \nu(dz) + \int_{[-1,1]^c}z\nu(dz).
 \eeqnn

Clearly, the two L\'evy processes $\{L(t): t\ge 0\}$ and $\{\xi(t): t\ge 0\}$ generate the same filtration. Let $\mbf{P}^\xi$ denote the conditional law given $\{L(t): t\ge 0\}$ or $\{\xi(t): t\ge 0\}$. Let $Z(t) = X(t)\exp\{-\xi(t)\}$. Following the arguments in Bansaye et al.\ (2013) or He et al.\ (2016) it is not hard to check that
 \beqlb\label{eq5.1}
\mbf{P}^\xi [e^{-\lambda Z(t)}|\mcr{F}_r] = \exp\{-Z(r)u^\xi_{r,t}(\lambda)\},
 \qquad
\lambda\ge 0,~t\ge r\ge 0,
 \eeqlb
where $r\mapsto u^\xi_{r,t}(\lambda)$ is the solution to
 \beqnn
\frac{d}{dr}u^\xi_{r,t}(\lambda) = ce^{-\alpha \xi(r)} u^\xi_{r,t}(\lambda)^{1+\alpha},
 \qquad
u^\xi_{t,t}(\lambda) = \lambda.
 \eeqnn
By solving the above equation, we get
 \beqlb\label{eq5.2}
u^\xi_{r,t}(\lambda)
 =
\Big(c\alpha\int_r^t e^{-\alpha \xi(s)}ds + \lambda^{-\alpha}\Big)^{-1/\alpha}.
 \eeqlb
From (\ref{eq5.1}) and (\ref{eq5.2}) we see that the survival probability of the CBRE-process up to time $t\ge 0$ is give by
 \beqlb\label{eq5.3}
\mbf{P}(X(t)>0)
 =
\lim_{\lambda\to \infty}\mbf{P}\big[1-\exp\{-xu^\xi_{0,t}(\lambda)\}\big]
 =
\mbf{P}\Big[F_x\Big(\int_0^t e^{-\alpha \xi(s)}ds\Big)\Big],
 \eeqlb
where $F_x(z) = 1 - \exp\{-x(c\alpha z)^{-1/\alpha}\}$. Let ${\it\Phi}(\lambda) = \log \mbf{P}\exp\{\lambda \xi(1)\}$ denote the Laplace exponent of $\{\xi(t): t\ge 0\}$.

The following theorem is an immediate consequence of Theorem~\ref{th2.9}. Using the notation introduced there, it gives characterizations of the five regimes of the asymptotics of the survival probability of the CBRE-process:

\begin{theorem}\label{th5.1} Suppose that $\{0,1\}\subset \mcr{D}^\circ({\it\Phi})$. Then we have the following five regimes of the survival probability of the CBRE-process:
 \begin{enumerate}

\item[{\rm(1)}] {\rm(Supercritical case)} If $0<{\it\Phi}'(0)$, we have the non-degenerate limit
 \beqnn
\mbf{P}(\tau_0< \infty) = \lim_{t\to \infty} \mbf{P}(X(t)=0) = \mbf{P}[F_x(A_\infty^\alpha(\xi))].
 \eeqnn

\item[{\rm(2)}] {\rm (Critical case)} Suppose that Condition~\ref{th2.8} is satisfied and ${\it\Phi}'(0) = 0$. Then we have the non-degenerate limit
 \beqnn
\lim_{t\to \infty} t^{1/2}\mbf{P}(X(t)>0) = \sqrt{\frac{2}{\pi{\it\Phi}''(0)}}\hat{\mbf{P}}[H(1)] D_2(\alpha,F_x).
 \eeqnn

\item[{\rm(3)}] {\rm (Weakly subcritical case)} Suppose that Condition~\ref{th2.8} is satisfied and ${\it\Phi}'(0)< 0< {\it\Phi}'(1)$. Let $\varrho\in (0,1)$ be the solution of ${\it\Phi}'(\varrho)=0$. Then we have the non-degenerate limit
 \beqnn
\lim_{t\to \infty}t^{3/2}e^{-t{\it\Phi}(\varrho)}\mbf{P}(X(t)>0)
 =
\frac{c(\varrho)}{\sqrt{2\pi{\it\Phi}''(\varrho)}} D_3(\alpha,F_x).
 \eeqnn

\item[{\rm(4)}] {\rm (Intermediately subcritical case)} Suppose that Condition~\ref{th2.8} is satisfied and ${\it\Phi}'(1) = 0$. Then we have the non-degenerate limit
 \beqnn
\lim_{t\to \infty}t^{1/2}e^{-t{\it\Phi}(1)}\mbf{P}(X(t)>0)
 =
x(c\alpha)^{-1/\alpha}\sqrt{\frac{2}{\pi{\it\Phi}''(\beta)}} \mbf{P}^{(1)}[H(1)] D_4(\alpha,1).
 \eeqnn

\item[{\rm(5)}] {\rm (Strongly subcritical case)} If ${\it\Phi}'(0)< 0$, we have the non-degenerate limit
 \beqnn
\lim_{t\to \infty}e^{-t{\it\Phi}(1)}\mbf{P}(X(t)>0)
 =
x(c\alpha)^{-1/\alpha} \mbf{P}^{(1)}[A_\infty^\alpha(-\xi)^{-1/\alpha}].
 \eeqnn

 \end{enumerate}
\end{theorem}

The convergence rate of the survival probability of CBRE-processes has been studied by B\"oinghoff and Hutzenthaler (2012), Bansaye et al.\ (2013), Palau and Pardo (2015a) and Palau et al.\ (2016). But, the limiting constant has only been computed explicitly in B\"oinghoff and Hutzenthaler (2012) for the case where the environment process is a Brownian motion with drift.

\end{document}